\documentclass[11pt,a4paper]{article}

\usepackage[english]{babel}
\usepackage{amsmath,amsthm,amssymb}
\usepackage{graphicx}

{
\newtheorem{theorem}{Theorem}

\newtheorem{lemma}{Lemma}
\newtheorem{conjecture}{Conjecture}
}
\newcommand{\beq}{\begin{equation}}
\newcommand{\eeq}{\end{equation}}
\newcommand{\beqs}{\begin{equation*}}
\newcommand{\eeqs}{\end{equation*}}

\newcommand{\reals}{\mathbb{R}}
\newcommand{\calDm}{\ensuremath{\mathcal{D}_{m}}}
\newcommand{\calL}{\ensuremath{\mathcal{L}}}
\newcommand{\calM}{\ensuremath{\mathcal{M}}}

\newcommand{\prob}[1]{\ensuremath{{\mathbb P}\!\left( #1 \right)}}
\newcommand{\expec}[1]{\ensuremath{{\mathbb E}\mspace{-1mu}\left[#1\right]}}
\newcommand{\pmu}[1]{\ensuremath{{\mathbb P}_{\mspace{-4mu}\mu}\!\left( #1 \right)}}
\newcommand{\pnu}[1]{\ensuremath{{\mathbb P}_{\mspace{-4mu}\nu}\!\left( #1 \right)}}
\newcommand{\ppi}[1]{\ensuremath{{\mathbb P}_{\mspace{-4mu}\pi}\!\left( #1 \right)}}
\newcommand{\psub}[2]{\ensuremath{{\mathbb P}_{\mspace{-4mu}#1}\!\left( #2 \right)}}
\newcommand{\supp}[1]{\ensuremath{{\rm supp}\left( #1 \right)}}

\newcommand{\md}{\mathrm{d}}

\newcommand{\norm}[1]{\ensuremath{\|{#1}\|}}

\newcommand{\ceil}[1]{\ensuremath{\lceil #1 \rceil}}

\input prepictex
\input pictexwd
\input postpictex

\title{Extremal distributions for tail probabilities of sums of iid random variables on $[0,1]$
\footnote{A preliminary version of this paper was presented at IWAP2008, Compi\`egne, France.}}
\author{Ludolf E.~Meester}
\date{\today}
\begin{document}
\maketitle
\begin{abstract}
\noindent
Two old conjectures from problem sections, one of which from \emph{SIAM Review}, concern the
question of finding distributions that maximize \prob{S_n\leq t}, where $S_n$ is the sum of
i.i.d.~random variables $X_1, \ldots, X_n$ on the interval $[0,1]$, satisfying $\expec{X_1}=m$.
In this paper a Lagrange multiplier technique is applied to this problem, yielding
necessary conditions for distributions to be extremal, for arbitrary $n$.
For $n=2$, a complete solution is derived from them: extremal distributions 
are discrete and have one of the following supports, depending on $m$ and $t$:
$\{0,t\}$, $\{t-1,1\}$, $\{t/2,1\}$, or $\{0,t,1\}$.
These results suffice to refute both conjectures. However, acquired insight naturally leads to a revised conjecture:
that extremal distributions always have at most three support points 
and belong to a (for each $n$, specified) finite collection of two and three point distributions.
\end{abstract}

\textbf{Keywords:} Probability theory, 
sums of iid random variables, 
Hoeffding inequality, 
extremal distributions.

\textbf{AMS Subject Classification:} 60E15.


\section{Two unsolved problems}
\label{twoproblems}
The problem section of the June 1986 issue of \emph{SIAM Review} lists the following, 
labeled \emph{Problem 86-6$^{*}$}~\cite{SRproblem866}:
\begin{quote}
In many audit populations items may have partial errors. Suppose each item in the population has an error size known
to be in the interval $[0,1]$. Suppose the mean population error is $m$ where $0<m<1$. A simple random sample of size
$n$ is drawn with replacement from that population. Let $S_n$ be the random variable representing the sum of the error
sizes of the $n$ sampled items. Given a constant $t <mn$ how should the error sizes be distributed in the population
to maximize \prob{S_n \leq t}?
It is conjectured that for each $m$ and $t$, there is a population with just two error sizes, one of which is $0$ or
$1$, such that \prob{S_n \leq t} is maximized. Prove or disprove.
\end{quote}

It is added that if this conjecture is true, then it will be possible to determine simple bounds on upper confidence
limits for some audit sampling problems.

The problem section of \emph{Statistica Neerlandica}, Vol.~47, no.~1,  lists the following as 
\emph{Problem 294}~\cite{SNproblem294}:
\begin{quote}
Consider i.i.d.~random variables $X_1, \ldots, X_n$ with $0 \leq X_i \leq 1$ and $\expec{X_i}=m$ given. 
Let $S_n=X_1+ \cdots +X_n$.  Consider the following statement:
\[ p= \prob{S_n \leq t} \quad \text{is maximal if} \quad \prob{X_1=1}=1-\prob{X_1=0}=m. \]
Show that this statement holds for all $t$ such that $p \leq p_0$ for some $p_0<1$ and find such a value for $p_0$.
\end{quote}
 
It appears that no solutions to these problems have been published.
This paper addresses them and presents a (partial) solution by considering:
\begin{quote}
Let $0<m<1$ and let \calDm\ be the set of probability measures on $[0,1]$ with mean $m$. Let $X_1, \ldots, X_n$ be
i.i.d.\ $\mu \in \calDm$ and $S_n=X_1+\cdots+X_n$. Determine $p(m,t) = \sup_{\mu \in \calDm} 
\pmu{S_n \leq t}$ and, if possible, (all) $\mu$ attaining the maximum.
\end{quote}
Note 
that $p(m,t)=1$ for $mn \leq t$: set $X_i = m$, all $i$;
for $mn \leq t$ the question of maximizing \prob{S_n \geq t} would be more natural. 
However, since $mn < t \leq n$ implies $n-t<n(1-m)$ and 
$\prob{S_n \geq t} = \prob{n-S_n \leq n-t} \leq p(1-m,n-t)$, this case is 
included in the problem statement; $t=mn$ once again is the trivial case.
So henceforth, $0 \leq t < mn$ is assumed.

In the sequel, when emphasis on the dependence on~$\mu$ is required,
$\pmu{\cdot}$ will be used. 
The supremum $p(m,t)$ is indeed attained by an element of \calDm, by Weierstra\ss' theorem, 
because the set \calDm\ is weak*-compact and $\mu \mapsto \pmu{S_n \leq t}$ is weak*-continuous.
A $\mu \in \calDm$ is called \emph{extremal} (for certain $n$, $m$ and $t$) if $\pmu{S_n \leq t}=p(m,t)$.

\medskip

The problem at hand satisfies a common rule: 
$n=1$ is trivial, $n=2$ can be solved with a reasonable amount of work, and $n \geq 3$ is hard.
After the $n=1$ case, we start with some general observations.
After that, some relevant results from the literature are discussed, which show that part of the $n=2$ case
follows from a paper by Hoeffding and Shrikande~\cite{hoeffdingshrikande} from~1955. We embark on a different
approach, applying a Lagrange multiplier technique from Mattner~\cite{mattner}, in Section~\ref{mattnersapproach},
for arbitrary $n$. 
The resulting Lagrange conditions provide a characterization of extremal distributions.
For the $n=2$ case, this allowed us to show that supports of extremal distributions necessarily look like 
$\{s, t-s\}$ or $\{s, t-s, 1\}$, for some $s$.
After this reduction, shown in Section~\ref{threesupportpoints}, the search for extremal distributions may be
restricted to this more manageable class and a complete analysis is carried out.
As it turns out, the conjectures stated in the problems above
can be refuted based on these results (Section~\ref{conjectures}). 
The first conjectured solution, however, seems ``almost true'' and its exception is understandable, 
so in Section~\ref{revised} a revised and sharpened conjecture is formulated, including the specification, for each
$n$, of a collection of distributions of which the extremal one is conjectured to be a member.

\subsection{The case $n=1$}
Markov's inequality implies that
\beq
\label{markov}
\prob{X \leq t} = \prob{1-X \geq 1-t} \leq \frac{1 - \expec{X}}{1-t} = \frac{1-m}{1-t} 
\eeq
and this upper bound is attained by the following two-point distribution:
$\prob{X=t}=(1-m)/(1-t)$ and $\prob{X=1}=(m-t)/(1-t)$.

\subsection{Some results from the literature}
Hoeffding and Shrikande~\cite{hoeffdingshrikande} obtained results on
the supremum of the distribution function of the
i.i.d.~sum of \emph{two} random variables, given $k$ moment conditions and a restricted range. 
They showed that the supremum over all such distributions is the same as that over all discrete distributions with
at most $2k+2$ support points. In addition, they provide the following bound for \emph{nonnegative} i.i.d.~$X_1$
and $X_2$ with $\expec{X_1}=\gamma$:
\beq
\label{hsbound}
\prob{X_1+X_2 \geq c \gamma} \leq 
\begin{cases}
1 & \text{if $c \leq 2$;}\\
4/c^2 & \text{if $2 \leq c \leq \frac52$;} \\
2/c- 1/c^2 & \text{if $\frac{5}{2} \leq c$.}
\end{cases}
\eeq
For i.i.d.~$X_1$ and $X_2$ on $[0,1]$, with $\expec{X_1}=m$, 
one may translate the above result to one on the \emph{left} tail,
by switching to the complements with respect to 1: 
\beqs
\prob{X_1+X_2 \leq t}
= \prob{{1-X_1}+{1-X_2} \geq 2-t},
\eeqs
and (\ref{hsbound}) applies with $\gamma= 1-m$ and $c=(2-t)/(1-m)$.
The distributions attaining the resulting bound have as their support, respectively, $\{m\}$, $\{t/2,1\}$, and $\{t-1,1\}$.
The first applies to $c\leq 2$, or $t\geq 2m$, the trivial case; 
the second to $2 \leq c \leq 5/2$, or $t/2 \leq m \leq (2t+1)/5$;
the third to $5/2 \leq c$, or: $t \geq 1$ and $(2t+1)/5 \leq m \leq 1$.
These results resolve the $n=2$ case for a \emph{subset} of $(m,t)$-values. 
Furthermore, Hoeffding and Shrikande~\cite{hoeffdingshrikande} 
did not address the question of uniqueness of the extremal distributions.

\medskip

Hoeffding's well known inequality (see \cite{hoeffding}) bounds
deviations from the expected value for the average of independent, not-necessarily identically distributed, random
variables. The author states that the bound is not optimal, but the best bound that can be obtained via his
method, based on the moment generating function. Applying Theorem~1~\cite{hoeffding}, one obtains:
\beq
\label{ineqHoeffding}
\prob{S_n \leq t} \leq \left ( \frac{1-m}{1-t/n} \right )^{n-t} \, \left (\frac{m}{t/n} \right)^t.
\eeq
For $n=1$ it is clear that (\ref{markov}) is sharper, since $m>t$.

\subsection{Some general observations}
\label{generalobs}
\paragraph{One may assume $0$ or $1$ to be in the support.}\label{zeroandone}
Recall the definition of the \emph{support} of a measure $\mu$: the smallest closed set with measure 1; we denote it
by \supp{\mu}.
Suppose, for some extremal $\mu$ one has $\supp{\mu} \subset [a,b]$ with $0<a<b<1$.
Let $X_i$ have distribution $\mu$ and $Z_i = m + \alpha (X_i -m)$, $i=1,\ldots, n$.
Then $\expec{Z_i}=m$ and there exist $\alpha>1$ such that $\supp{Z_i} \subset [0,1]$. 
Writing $S^*_n= Z_1+\cdots+Z_n=\alpha S_n-(\alpha-1)\,m n$ one has for $t <m n$:
\[ \prob{S^*_n \leq t}=\prob{S_n \leq \tfrac1\alpha\, t + \tfrac{\alpha-1}{\alpha} \,m n} \geq \prob{S_n \leq t}. \]
Thus, if $\mu$ is extremal, a measure $\mu^*$ can be found that is extremal as well, and if the largest $\alpha$ is
chosen that satisfies $\supp{\mu^*} \subset [0,1]$, then \supp{\mu^*} will contain $0$ or $1$.
\paragraph{The supremum $p(m,t)$ is non-increasing in $m$.}
Fix $n$ and $t$.  Suppose $m_1 < m_2$ and $\mu_2$ attains the maximum value $p(m_2,t)$. 
Suppose $\mu_2$ has distribution function $G$. Define, for $0 \leq r \leq 1$, the distribution
function $F_r$ by $F_r(x)=\max(r, G(x))$. 
As $r$ goes from $0$ to $1$ the expectation of the corresponding distribution decreases from $m_2$ to $0$,
continuously, so for some $r$ the corresponding distribution $\mu_1$ has expectation~$m_1$.
This measure is stochastically smaller than $\mu_2$.
Hence, $p(m_1,t) \geq \psub{\mu_1}{S_n \leq t} \geq \psub{\mu_2}{S_n \leq t} = p(m_2,t)$.
\paragraph{Subprobability measures with expectation at least $m$.}
Instead of taking the supremum over \calDm\ one could take the supremum over 
the set of all \emph{sub}probability measures on $[0,1]$ with mean \emph{at least} $m$,
and the same $p(m,t)$ would result. 
In order to show this, suppose $\nu$ is a subprobability measure with defect $\rho=1-\nu([0,1]) \geq 0$ and
$\int_0^1 x \, \nu(\md x) = r \geq m$. Then
\[
\pnu{S_n \leq t} \leq \psub{\nu+ \rho \cdot \delta_0}{S_n \leq t} \leq p(r,t) \leq p(m,t),
\]
where the first inequality states that putting the mass-defect in $0$ will lead to an improvement, and the last
inequality follows from the non-increasing\-ness proved above. Considering that $d>0$ or $r>m$ will make at least one
of the inequalities strict, it is clear that $p(m,t)$ can only be attained with $d=0$ and $r=m$. 
\paragraph{An upper confidence bound on $m$.} 
Using $S_n$ as test statistic one may define a non-parametric confidence bound on $m$, as follows.
Let $m_u$ be solution of $p_n(m,t)=\alpha$. By the non-increasing\-ness proved above,
this implies that $p_n(m,t) \leq \alpha$ for $m \geq m_u$. So, for 
$m \geq m_u$ it follows that $\prob{S_n \leq t} \leq p_n(m,t) \leq p_n(m_u,t)=\alpha$.

\section{Mattner's Lagrange approach}
\label{mattnersapproach}
Mattner~\cite{mattner} developed a general method for treating extremal problems for probability distributions. 
His main theorem is stated below and subsequently applied to the problem:
\begin{theorem} 
Let $Z$ be a Banach space, 
$\varphi_i : Z \rightarrow \reals$, $i=0, \ldots, k$ and
$\psi_j : Z \rightarrow \reals$, $j=1, \ldots, l$, continuously Fr\'echet-differentiable, 
and $C$ a convex cone in~$Z$.
Define the Lagrange functional
\[ \calL(z) := \lambda_0 \varphi_0(z)+ \sum_{i=1}^k \lambda_i \varphi_i(z) + \sum_{j=1}^l \alpha_j \psi_j(z), \]
and let $\partial\calL(z;w)$ denote the Fr\'echet-derivative of $\calL(z)$ in direction $w$.
If $z \in Z$ minimizes $\varphi_0$ subject to
\begin{align*}
\varphi_i(z) &\leq 0, \quad i=1, \ldots, k,\\
\psi_j(z) &= 0,   \quad j=1, \ldots, l,\\
z & \in C, 
\end{align*}
then there exist $\lambda_0, \lambda_1, \ldots, \lambda_k, \alpha_1, \ldots, \alpha_l \in \reals$ with
\begin{tabbing}
(i)\quad    \= not all $\lambda_i$ and $\alpha_j$ vanish,\\
(ii)   \> $\lambda_i \geq 0$, $i=0, \ldots, k$, \\
(iii)  \> $\partial\calL(z;w) \geq 0$, $w \in C$ \\
(iv)   \> $\partial\calL(z;z) = 0$.
\end{tabbing}
\end{theorem}
\paragraph{Application to the problem.} 
Let \calM\ be the set of signed Borel measures on $[0,1]$, with
norm $\norm{\mu}=\int |\mu(\md x) |$ (total variation). The pair $(\calM, \norm{\,\cdot\,})$ is a
Banach space and probability measures are contained in the positive 
cone~$C=\{\mu \in \calM: \mu(B) \geq 0 \: \text{for every Borel set $B$}\}$. 
Define
\beqs
\varphi_0(\mu) = -\prob{S_n \leq t}= - \int_A \mu(\md x_1) \cdots \mu(\md x_n),
\eeqs
where $A=\{(x_1,x_2, \ldots, x_n): x_1+\cdots+x_n \leq t\}$. 
The (clearly continuous) Fr\'echet-derivative of $\varphi_0$ is given by
\beqs
\partial\varphi_0(\mu;\nu)=-n \, \int_A  \mu(\md x_1) \cdots \mu(\md x_{n-1}) \nu(\md x_n) =
-n \int_0^1 \prob{S_{n-1}\leq t-x} \, \nu(\md x).
\eeqs
Indeed, this follows from the next formula, which is established by binomial expansion of the $n$-fold 
product of the measure $\mu+\nu$:
\beqs
\norm{\varphi_0(\mu+\nu)-\varphi_0(\mu)-\partial\varphi_0(\mu;\nu)} = O\left(\norm{\nu}^2\right).
\eeqs
For the constraints define
\beqs
\varphi_1(\mu) = \mu([0,1]) -1\quad\text{and}\quad \varphi_2(\mu) =m - \int_0^1 x \, \mu(\md x),
\eeqs
whose (continuous) Fr\'echet-derivatives are given by
\beqs
\partial\varphi_1(\mu;\nu) =\nu([0,1]) \quad\text{and}\quad \partial\varphi_2(\mu;\nu) = - \int_0^1 x \, \nu(\md x).
\eeqs
Now, define the Lagrange functional:
$\calL(\mu) = \lambda_0\, \varphi_0(\mu) + \lambda_1\, \varphi_1(\mu) + \lambda_2\, \varphi_2(\mu)$.

From Mattner's theorem one concludes: 
if $\mu$ minimizes $\varphi_0$ subject to $\varphi_1(\mu) \leq 0$,
$\varphi_2(\mu) \leq 0$, $\mu \geq 0$, then there exist nonnegative $\lambda_0$, $\lambda_1$, and $\lambda_2$, not all
zero, such that
\begin{align}
\label{Lcond1}
\partial\calL(\mu;\nu) & \geq 0, \quad \text{for $\nu \geq 0$,}\\
\label{Lcond2}
\partial\calL(\mu;\mu) & =0,
\end{align}
where $\partial\calL(\mu;\nu)$ is the Fr\'echet-derivative of $\calL$ at $\mu$ in direction $\nu$
and given by
\beqs
\partial\calL(\mu;\nu) = \int_0^1 \ell(x) \, \nu(\md x)
\eeqs
with
\beqs
\ell(x)= -n \, \lambda_0 \, \prob{S_{n-1} \leq t-x} + \lambda_1 -\lambda_2 \, x.
\eeqs
Note that $\ell(x)$ is continuous from the left and that jump-discontinuities (if any) are upwards.

\subsection{The Lagrange conditions}
\label{Lagrange}
From Mattner's theorem some properties of extremal distributions can be derived, as well as an expression of
\prob{S_n\leq t} in terms of the Lagrange multipliers. First, the redundant Lagrange multiplier
$\lambda_0$ is removed.

From Lagrange condition (\ref{Lcond1}), by substituting $\nu = \delta_x$ (point-mass at $x$), one may conclude
$\ell(x) \geq 0$, for $0 \leq x \leq 1$. Combining this with the second Lagrange condition (\ref{Lcond2}) results in:
\beq
\label{lzeroae}
\ell(x)=0 \quad \text{for $\mu$-a.e.~$x$}.
\eeq
It is first argued that $\lambda_0$ cannot be zero. 
If $\lambda_0=0$, then $\ell(x)=\lambda_1 -\lambda_2 \, x$ should be nonnegative for 
$0 \leq x \leq 1$, whence $\lambda_1 \geq \lambda_2 \geq 0$ and, necessarily, $\lambda_1 >0$, for they cannot all
three be zero. However, $\ell(x)=0$ must
have at least one solution, or else $\supp{\mu}=\emptyset$. This leaves $\lambda_1=\lambda_2$ as sole possibility,
implying that $\mu=\delta_1$, which contradicts the assumption $\expec{X}=m <1$.
Therefore, $\lambda_0>0$ and without loss of generality it is henceforth assumed that $n \lambda_0=1$.

Lagrange condition~(\ref{Lcond2}), $\varphi_1(\mu)\leq 0$, and $\varphi_2(\mu)\leq 0$, imply
\beq
\label{probbound}
\prob{S_n\leq t} = \lambda_1 \mu([0,1])-\lambda_2 \int_0^1 x \, \mu(\md x).
\eeq
The Lagrange conditions can be restated as
\begin{align}
\label{L1}
\prob{S_{n-1} \leq t-x } & \leq \lambda_1 -\lambda_2 \, x, \quad 0 \leq x\leq 1,\quad
\text{and}\\
\label{L2}
\prob{S_{n-1} \leq t-x } & = \lambda_1 -\lambda_2 \, x, \quad\text{for $x \in \supp{\mu}$}.
\end{align}
The following lemma shows
that the last statement follows from (\ref{lzeroae}):
\begin{lemma} 
\label{elliszero}
Let $\mu$ be extremal. Then $\ell(x)=0$ for $x \in \supp{\mu}$.
\end{lemma}
\textbf{Proof.}
Let $x \in \supp{\mu}$ and suppose a Borel set $A \subset [0,1]$ satisfies $\mu(A)=1$ and $\ell(y)=0$ for $y \in A$.
If $x$ is an atom of $\mu$ then $x \in A$ and $\ell(x)=0$ follows.
Otherwise, if $x$ is an interior point or a right boundary point 
of \supp{\mu}, a sequence $(x_k)$ can be found within $A$ such that $x_k \uparrow x$,
whence $\ell(x)=0$, by left-continuity of $\ell$. If $x$ is a left boundary point, 
one can find within $A$ a sequence $x_k \downarrow x$, whence $0\leq \ell(x)\leq \ell(x+)=\lim
\ell(x_k)=0$, since jumps cannot go down.~\qed

\medskip

It is shown that $\lambda_2>0$ must hold.
Let $s=\min \supp{\mu}$ and $u=\max \supp{\mu}$, 
then $s \leq m \leq u$ and gaps in \supp{S_{n-1}} cannot exceed $u-s$ in length.
Lagrange conditions (\ref{L1}) and (\ref{L2}) imply $\prob{t-u < S_{n-1} \leq t} \leq \lambda_2 \, u$.
The probability, however, must be positive:
$\expec{S_{n-1}}= (n-1)\,m > t -m \geq t-u$ implies that \prob{t-u < S_{n-1}} must be positive;
$\prob{t < S_{n-1}}=1$ cannot be the case, or else $\prob{S_{n}\leq t}=0$ and $\mu$ is not extremal.
\paragraph{Support conditions.} The Lagrange conditions imply several properties for the support of $S_{n-1}$
and $S_{n}$.  An immediate consequence of the next lemma is that $t \in \supp{S_n}$.
\begin{lemma} 
\label{support}
Let $\mu$ be extremal, $x \neq 1$.
If $x \in \supp{\mu}$ then $t-x \in \supp{S_{n-1}}$.
\end{lemma}
\textbf{Proof.} By contraposition. 
Suppose $t-x \not\in \supp{S_{n-1}}$, for some $0 < x < 1$. 
Let $B_\epsilon(x)=(x-\epsilon,x+\epsilon)$.
Then $\prob{S_{n-1} \in B_\epsilon(t-x)}=0$ for some $\epsilon>0$, 
implying that $\prob{S_{n-1}\leq t-y}$ is constant for $y \in B_\epsilon(x)$ 
and that $\ell(y) \geq 0$ is linearly decreasing on this set, 
which implies $\ell(x)> 0$. 
For $x=0$, this reasoning shows that $\ell(y)$ is linearly decreasing for $y \in [0,\epsilon)$, 
with $\ell(0)>0$ as conclusion.
For $x=1$, nothing about the positivity of $\ell(1)$ can be concluded from the fact that $\ell(y)$ 
is linearly decreasing and positive for $y \in (1-\epsilon,1]$; $\ell(1)=0$ is still possible.
So, for $0\leq x <1$, $t-x \not\in \supp{S_{n-1}}$ implies $\ell(x)>0$, 
which by Lemma~\ref{elliszero} implies $x \not \in \supp{\mu}$.~\qed

\medskip

The Lagrange conditions in this section provide necessary conditions (\ref{L1}), (\ref{L2}), and
Lemma~\ref{support}, that should be satisfied by extremal distributions. 
It is not difficult, for general $n$, to identify a number of distributions that satisfy them
(see Section~\ref{revised}). 
However, unless \emph{all} the solutions are identified, there are no guarantees that the best of
the solutions found indeed attains the supremum $p(m,t)$.

\section{The case $n=2$}
\label{threesupportpoints}
Lemma~\ref{support} yields an especially strong result for $n=2$, because $S_{n-1}=X_1$ and the lemma
characterizes the support of (candidate) extremal distributions.
Below, certain two and three point solutions to the Lagrange conditions will be identified.
Other solutions (if any) cannot be extremal:
it will be shown that one can always find a distribution of the two or three point type that has a strictly
larger \prob{S_n \leq t}-value.  Hence, all extremal distributions belong to this special class.

\medskip

Suppose $\mu$ satisfies the Lagrange conditions and $s = \min \supp{\mu}$. 
Lemma~\ref{support} implies that $t-s \in \supp{\mu}$, and if this is not the largest support point, 
then $\max \supp{\mu}=1$. Therefore, two cases are to be considered.


\medskip

First, assume that $t-s=\max \supp{\mu}$. Note that, necessarily, $0 \leq s \leq t-s \leq 1$ and $t-s
>m$ (or else  $\int_0^1 x \mu (\md x) < m$), which imply $s <t/2 <t-s$ by $m >t/2$. 
Note that $m<t$ must hold, or no such $\mu$ exist.

Let $F$ be the distribution function corresponding to $\mu$.
Lagrange condition~(\ref{L2}) requires that nonnegative $\lambda_1$ and $\lambda_2$ exist such that:
\beqs
F(t-s) = \lambda_1 -\lambda_2 \, s \quad \text{and} \quad 
F(s) = \lambda_1 -\lambda_2 \, (t-s).
\eeqs
From the monotonicity and nonnegativity of $F$:
\beq
\label{intFx}
\int_0^1 F(x) \md x \geq (t-2s) \, F(s) + (1-t+s) \, F(t-s),
\eeq
where equality holds (if and) only if $s$ and $t-s$ are the only support points. 
Since $F(t-s)=F(1)$, one may write $\lambda_1=F(1)+\lambda_2 \, s$ and $F(s)=F(1)-\lambda_2 \, (t-2s)$.
Combining things, one obtains:
\beq
\int_0^1 x \mu \, (\md x) = F(1) - \int_0^1 F(x) \, \md x 
\leq s \, F(1) + \lambda_2 \, (t-2s)^2
\eeq
whence $\lambda_2 \geq (m-s)/(t-2s)^2$.
Starting from~(\ref{probbound}), this results in:
\beq
\label{Pbound2}
\pmu{S_2 \leq t} \leq \lambda_1  - \lambda_2 \, m = F(1) - \lambda_2 (m-s) 
\leq 1 - \left ( \frac{m-s}{t-2s} \right )^2.
\eeq
Note that $0 < (m-s)/(t-2s) < 1$ since $s <t/2<m$ and $m <t-s$. Let~$\pi$ be the probability 
measure $\pi$ on $\{s,t-s\}$ defined by $\pi_{t-s}= (m-s)/(t-2s)=1-\pi_s$, where $\pi_x:=\pi(\{x\})$. 
Then $\int_0^1 x \, \pi(\md x)=m$ and $\ppi{S_2\leq t}$ equals the right
hand side of~(\ref{Pbound2}). 
The upper bound on \pmu{S_2 \leq t} is strict, unless $F(1)=1$, $\int_0^1 x \, \mu(\md x)=m$, and equality
holds in~(\ref{intFx}). These conditions, however, uniquely identify $\pi$, showing that $\mu$ can only be extremal
if $\mu=\pi$.

\bigskip

Next, consider the situation where $t-s<\max \supp{\mu}=1$. Necessarily, $0 \leq s \leq t-s < 1$ must hold. 
Lagrange condition~(\ref{L2}) specifies for the support points $s$, $t-s$ and $1$, respectively:
\beqs
F(t-s) = \lambda_1 -\lambda_2 \, s, \quad
F(s) = \lambda_1 -\lambda_2 \, (t-s), \quad \text{and} \quad 
F(t-1) = \lambda_1 -\lambda_2.
\eeqs
Since $t-1<s=\min \supp{\mu}$, $F(t-1)=0$ and so $\lambda_1=\lambda_2$, which is used to eliminate $\lambda_1$.

Note that $F(t-s) < F(1)$. Further, that $F(1)=1$ must hold, or
the defect could be added as an atom in $0$, which would strictly enlarge $\pmu{S_2 \leq t}$.
If $F(1)=1$ and $\int_0^1 x \mu(\md x) > m$ then a small mass $\epsilon>0$ could be moved from~$1$
to $0$, still keeping the mean above $m$. This would increase \pmu{S_2 \leq t} by at least $\epsilon^2$.
Hence, if $\mu$ is to be extremal, then $F(1)=1$ and $\int_0^1 x \mu(\md x) = m$ must hold.

Combining (\ref{intFx}) with $F(t-s) = \lambda_2 \, (1-s)$ and $F(s) = \lambda_2 \, (1-t+s)$, one obtains
\beq
\label{intFx2}
1-m = \int_0^1 F(x) \, \md x \geq \lambda_2 \, (1-t+s) (1+t-3s),
\eeq
where equality holds (if and) only if $s$, $t-s$ and $1$ are the only support points. Apparently,
\beqs
\lambda_2 \leq \lambda_2^+ := \frac{1-m}{(1-t+s) (1+t-3s)}.
\eeqs
Define the measure $\pi$ on $\{s,t-s,1\}$ by
\beq
\label{threepointsol}
\pi_s=\lambda_2^+ (1-t+s), \quad
\pi_{t-s}=\lambda_2^+ (t-2s), \quad
\pi_1=1-\lambda_2^+ (1-s).
\eeq
(Note that $s=t-s$ leaves a valid probability measure on the set $\{t/2,1\}$.)
If $\lambda_2^+ (1-s) <1$, then $\pi$ is a probability measure with mean $m$:
the probabilities are nonnegative and sum to 1, 
$s \, \pi_s + (t-s) \, \pi_{t-s} +\pi_1= 1 - \lambda_2^+ \, (1-t+s) (1+t-3s)=m$. 
Furthermore, $\pi$ satisfies the Lagrange conditions (for $\lambda_2^+$) and so
\beqs
\pmu{S_2 \leq t} =
\lambda_2 \, (1-m)  \leq 
\lambda_2^+ \, (1-m) =
\ppi{S_2 \leq t}  = \frac{(1-m)^2}{(1-t+s) (1+t-3s)}.
\eeqs
The inequality is strict unless $\mu=\pi$; this can be seen from~(\ref{intFx2}).

If $\lambda_2^+ (1-s) \geq 1$, then
\beqs
(1-m)(1-s) \geq (1-t+s) (1+t-3s) = (1-s)^2 -(t-2s)^2,
\eeqs
which is equivalent to
\beq
\label{twopointonly}
(t-2s)^2 \geq (1-s)(m-s).
\eeq
Since $\lambda_2 \, (1-s) <1$ by assumption, starting from~(\ref{probbound}) (with $\lambda_1=\lambda_2$),
\beq
\label{Pbound3}
\pmu{S_2 \leq t} = \lambda_2 \, (1-m) < \frac{1-m}{1-s} = 1- \frac{m-s}{1-s}  
\leq 1 - \left ( \frac{m-s}{t-2s} \right )^2,
\eeq
where the last inequality follows from~(\ref{twopointonly}). 
That inequality also implies $t-2s>m-s$, which combined with $m>s$ guarantees the existence of the
probability measure with support $\{s, t-s\}$ and mean $m$. As was shown before, this measure attains 
the value on the right hand side of (\ref{Pbound3}).
In all cases it has now been shown that if $\mu$ satisfies the Lagrange conditions it equals a discrete measure on
$\{s, t-s\}$ or $\{s,t-s,1\}$, for some $s$, or else $\pmu{S_2 \leq t} < \ppi{S-2 \leq t}$ for some $\pi$ from this
class.

The last steps of the solution consist of optimizing within the class of two and three point support distributions just
identified. 

\bigskip

\textbf{Remark: an alternative approach?}\label{alternative}
What follows is a sketch of a proof that would work if one could show that
extremal measures cannot have a singular component. The Lagrange conditions imply that if the support contains an
interval, say $A$, then $\mu$ has density equal to $\lambda_2$ on that interval. Lemma~\ref{support} implies that the
same holds for $t-A$, from which it easily follows that $\mu$ can be improved upon by moving the mass to the center
of the intervals. If $\mu$ is purely atomic, a similar argument that exploits the symmetry of the support can be used
to show that $[0,t/2)$ cannot contain more than one atom. After this, four possible support points remain: (some)
$s$, $t/2$, $t-s$, and $1$. A simple mass transfer argument shows that the first three cannot occur together. This
leaves one with the same possibilities as in the current line of reasoning.

\subsection{$\{s,t-s\}$-solutions}
Recall that the bound from~(\ref{Pbound2}) and (\ref{Pbound3}) can be attained by the probability
measure $\pi$ on $\{s,t-s\}$ defined by $\pi_{t-s}= (m-s)/(t-2s)=1-\pi_s$. 
The largest \ppi{S_2 \leq t}-value is attained  for the smallest feasible $s$, as
$\pi_{t-s}$ is increasing in $s$. For $t \leq 1$ the maximum is at $s=0$, for $t \geq 1$ at $s=t-1$. Thus, the best
solutions of this type are as follows.

For $m < t \leq 1$: $\pi_0=1-m/t$, $\pi_t=m/t$ and
\beqs
\label{case2a}
\ppi{S_2 \leq t} =  1 - \left ( \frac{m}{t} \right)^2.
\eeqs
For $t\geq 1$: $\pi_{t-1}=(1-m)/(2-t)$, $\pi_1=(1+m-t)/(2-t)$ and
\beqs
\label{case2b}
\ppi{S_2 \leq t} =  1 - \left ( \frac{1+m-t}{2-t} \right)^2.
\eeqs

\subsection{$\{s,t-s,1\}$-solutions}
Candidate extremal distributions are the probability measures with support $\{s,t-s,1\}$ and probabilities given 
by~(\ref{threepointsol}), provided $0 \leq s \leq t/2$, $s \geq t-1$ and $\lambda_2^+ (1-s) \leq 1$, or:
\beq
\label{existence3}
(1-m) (1-s) \leq (1-t+s) (1+t -3s).
\eeq
For the sake of a simpler exposition two small additions were made to the class considered:
equality in the previous formula corresponds to boundary cases with $\pi_1=0$, which, just as the case $s=t-1$ that
was added, leads to distributions already considered. 

Recall that $\ppi{S_2 \leq t} = (1-m)^2/(1-t+s)(1+t-3s)$ for $\pi$ as in~(\ref{threepointsol}). What
remains is maximize over feasible $s$.
Define $p(s)=(1-t+s) (1+t -3s)$. Since $p$ is a concave function, the solutions
to~(\ref{existence3}) constitute an interval; call the left end point $s_0$. Since $p(s)=(1-s)^2
-(t-2s)^2$, (\ref{existence3}) is equivalent to $(t-2s)^2 \leq (1-s)(m-s)$, which shows that $s=t/2$ is always
feasible, and only feasible $s$ between $s_0$ and $t/2$ need to be considered.

Consider the maximization problem: since $p$ is concave, the maximum of $(1-m)^2/p(s)$ is attained at an
end point of the feasible range, i.e., $t/2$ or the left end point. Note, however, that as $s \downarrow s_0$, also $\pi_1
\downarrow 0$, and what results is a distribution on $\{s_0, t-s_0\}$, already considered in the previous section. 
This means that if $s_0 \geq 0$ and $s_0 \geq t-1$, the entire range $s_0 \leq s \leq t/2$ corresponds to feasible
solutions, at the left end dominated by solutions already considered. Then, the distribution corresponding
to the right end point is the only new (candidate) extremal distribution. It is given by: $\pi_{t/2}=(1-m)/(1-t/2)$ and
$\pi_1=(m-t/2)/(1-t/2)$ with
\beqs
\label{case2c}
\ppi{S_2 \leq t} =   \left ( \frac{1-m}{1-t/2} \right)^2.
\eeqs
Note that this solution exists for all $(m,t)$-pairs under consideration.
What remains now, is to determine whether the maximum can be attained for an intermediate value $s_0<s<t/2$, which
would correspond to a \emph{true} three point distribution.

\medskip

First, a closer look at $p(s)$ is warranted. It has zeros at $t-1$ and $(t+1)/3$, a maximum value of $(2-t)^2/3$ attained at
$(2t-1)/3$; $p(0)=1-t^2$ and $p((5t-4)/6)=p(t/2)=(1-t/2)^2$. These points are ordered in the following manner:
\beqs
t-1 < \frac{5t-4}6 < \frac{2t-1}3 <\frac{t}2 < \frac{t+1}3,
\eeqs
where zero can be anywhere to the left of $t/2$, depending on $t$.

\medskip

From equality in (\ref{existence3}) one sees that $p(s_0)=(1-m)(1-s_0)>0$, and since $t-1$ is the left zero of $p$,
this implies that $s_0>t-1$. Therefore, for $t\geq 1$, the whole range $s_0 \leq s \leq t/2$ corresponds to 
feasible solutions dominated by one of the two point solutions corresponding to the end points, i.e., 
with support $\{s_0,t-s_0\}$ or $\{t/2,1\}$. 
Next, consider $t <1$. A~true three point solution occurs if $s_0<0$, which happens if inequality
(\ref{existence3}) is strict 
for $s=0$, which happens if $t < \sqrt{m}$. However, if $t<4/5$, then $(5t-4)/6<0$ and $p(0)>p((5t-4)/6)=p(t/2)$, whence
$(1-m)^2/p(s)$ attains a higher value at $s=t/2$ than at $s=0$. In summary, this shows that the best solution is
obtained at $s=0$ only for $4/5 \leq t < \sqrt{m}$. It is the extremal probability measure $\pi$ on $\{0,t,1\}$ defined by:
$\pi_0 ={(1-m)(1-t)}/(1-t^2)$,  $\pi_t ={(1-m)\,t}/(1-t^2)$, and $\pi_1 =(m-t^2)/(1-t^2)$, with
\beq
\label{case3}
\ppi{S_2 \leq t} = \frac{(1-m)^2}{1-t^2}.
\eeq
The $(m,t)$-range where this distribution dominates all others has just been determined. 
On the complement of this range several two point solutions may exist together and therefore
need to be compared.

\subsection{Some comparisons}
In the region $m \leq t \leq 1$, both the $\{0,t\}$ and the $\{t/2,1\}$-solution exist. The second is the best when
\beq
\label{compareCandA}
1 - \left ( \frac{m}{t} \right )^2 < \left ( \frac{1-m}{1-t/2} \right )^2.
\eeq
Substituting $m=a t$, the equivalent inequality $4(1-at)^2 -(1-a^2)(2-t)^2 >0$ is obtained, 
which in turn simplifies to
$\left( 5\,{t}^{2}-4\,t+4 \right) {a}^{2}-8\,at-{t}^{2}+4\,t>0$.
The discriminant of this quadratic in $a$ is $20\,{t}^{4}-96\,{t}^{3}+144\,{t}^{2}-64\,t$ which factors as
$4\, t\, ( 5\,t-4 )\, (2-t) ^{2}$. This shows that the inequality~(\ref{compareCandA}) is valid for $0<t<4/5$,
as the discriminant is negative for these values.
For $4/5 \leq t \leq 1$, the boundary curve of the inequality~(\ref{compareCandA}) is given by
\beq
\label{monet}
m_1(t)= \frac {4\,{t}^{2}- \left( 2-t \right) t\sqrt {t \left( 5\,t-4
 \right) }}{5\,{t}^{2}-4\,t+4}.
\eeq
For $m<m_1(t)$ the $\{t/2,1\}$ solution is superior; for larger $m$ the $\{0,t\}$ solution is.

\medskip

To determine for which $m$ and $t$ the $\{t/2,1\}$-solution is best 
for $t \geq 1$ one needs to solve
\beqs
1-\left(\frac{m+1-t}{2-t}\right)^2 <  \left ( \frac{1-m}{1-t/2} \right )^2.
\eeqs
Setting $a=(1-m)/(2-t)$, this becomes $4a^2>1-(1-a)^2$, resulting in $a>2/5$, or $5m<1+2t$.

\bigskip

\begin{figure}
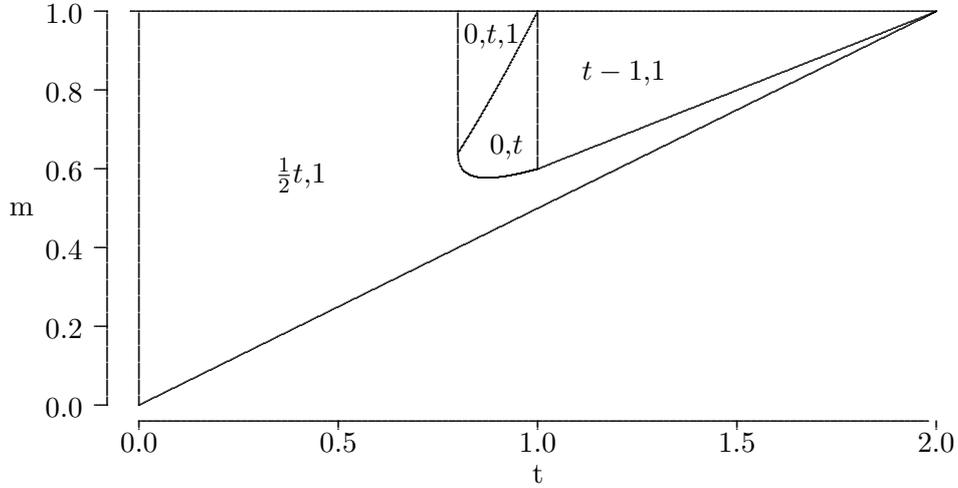

\input plotnistwo.pic
\caption{Support points of extremal distributions, for $n=2$.}
\label{supportpoints}
\end{figure}

Summarizing everything, one obtains the following table. 
The function $m_1(t)$ on the second line is given in equation~(\ref{monet}).
Figure~\ref{supportpoints} shows the regions with the support of the respective extremal distributions.

\begin{center}
\begin{tabular}{lcl}
       support & $(t,m)$-region & \prob{S_2 \leq t}\\[.5ex] \hline
$\{0, t\}$      & $\tfrac45 \leq t \leq 1,\, m_1(t) \leq m \leq t^2$   & $1 - \left ( \frac{m}{t} \right )^2$,
\\[.5ex]
$\{0, t, 1\}$ & $\tfrac45 \leq t < \sqrt{m}$ &  $\frac{(1-m)^2}{1-t^2}$, \\[.5ex]
$\{t-1, 1\}$    & $1 \leq t \leq 2$,\, $5m>2t+1$                     & $1-\left(\frac{m+1-t}{2-t}\right)^2$,
\\[.5ex]
$\{t/2, 1\}$    & everywhere else   &   $\left ( \frac{1-m}{1-t/2} \right )^2$.
\end{tabular}
\end{center}

\begin{figure}
\includegraphics[width=\textwidth]{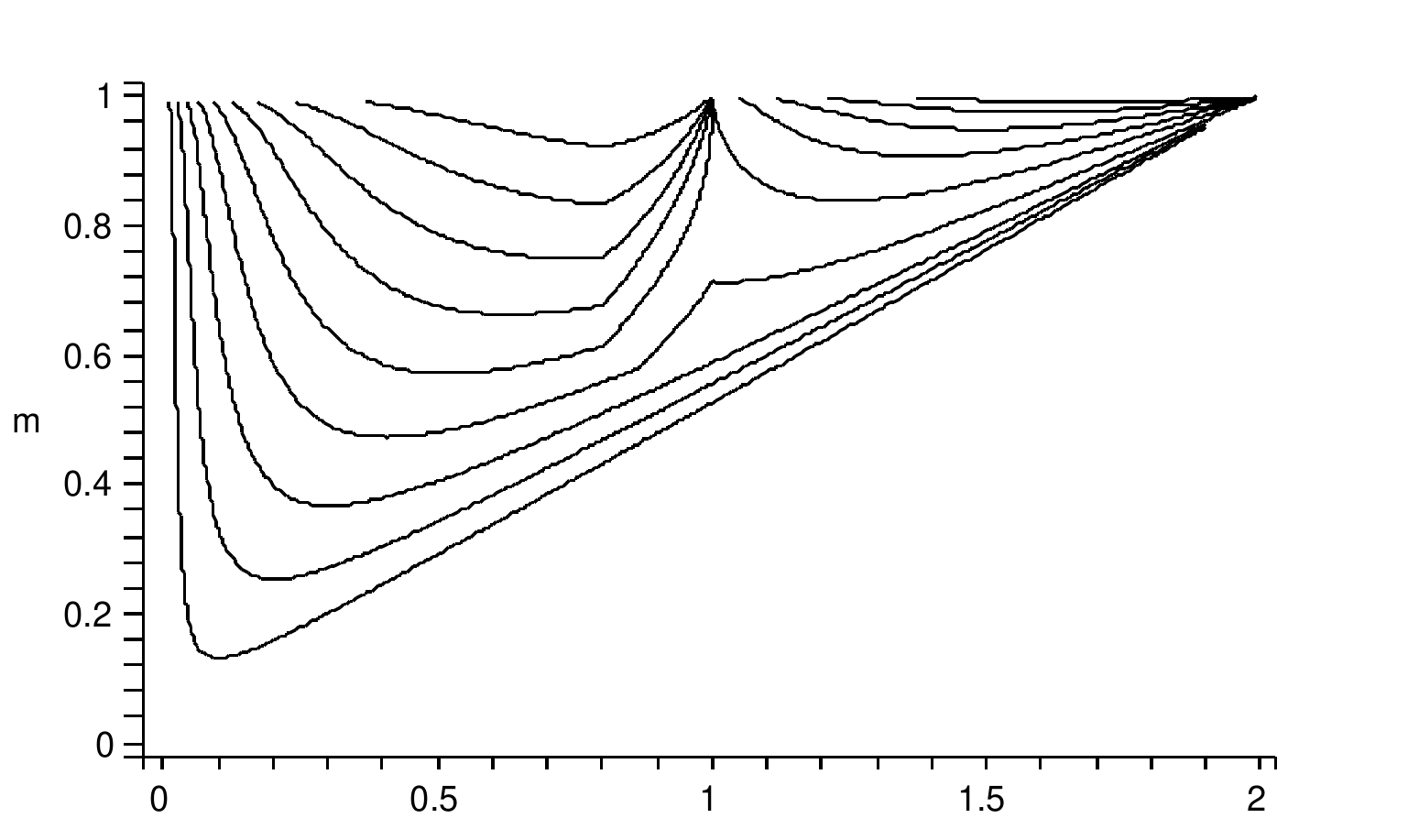}
\caption{Ratio of $p_2(m,t)$ and the Hoeffding bound~(\ref{ineqHoeffding}); 
contour lines correspond to $0.1$, $0.2$, \ldots, $0.9$, going from top to bottom.}
\label{versusHoeffding}
\end{figure}

The extremal measures are unique except on the boundary between the $\{t/2,1\}$-solution and the others, where two distinct
solutions yield the same \prob{S_2\leq t}-value; on the other boundaries, the two \emph{solutions} coincide.
Figure~\ref{versusHoeffding} shows a contourplot of the ratio of the supremum $p_2(m,t)$ and 
the Hoeffding bound~(\ref{ineqHoeffding}); the bound is sharp at the boundary $m=t/2$ and 
progressively looser as $m$ increases.

\section{The conjectures: one refuted, one revised}
\label{conjectures}
Strictly speaking, the results for $n=2$ suffice to disprove \emph{both} conjectures. 
Whereas the \emph{Statistica Neerlandica} conjecture can be utterly disproved,
the \emph{SIAM Review} conjecture is only disproved by a small $(m,t)$-region where
the extremal distribution has three support points (including $0$ and $1$!). 
In our view the \emph{SIAM Review} conjecture is close to what may be true and
therefore a revised conjecture is formulated below.

\paragraph{Statistica Neerlandica.}
It seems that what was meant is ``If $p_n(m,t) = \sup \prob{S_n \leq t}$ is small (enough), then the Bernoulli with
mean $m$ is the best distribution.'' 
Looking at the $n=2$ results, the Bernoulli only appears as extremal distribution for $t=0$ and
for $t=1$, $3/5 \leq m \leq 1$.
Let $b(n,p,x)$ denote the probability that a binomial random variable with parameters $n$ and $p$ attains a value less
than or equal to $x$. 
The following is the logical negation of the \emph{Statistica Neerlandica} conjecture:
\begin{lemma}
For any $0 < p_0 <1$ there exist $n$, $t$, and $m$, such that $b(n,m,t) < p_n(m,t) \leq p_0$.
\end{lemma}
\textbf{Proof.} Set $n=2$ and choose any $t$ in $(0,1)$ or $(1,2)$. 
Then for $m >t/2$: $b(n,m,t) < p_n(m,t)$ and as $m \uparrow 1$, $p_n(m,t) \to 0$.~\qed

\subsection{The \emph{SIAM Review} conjecture revised}
\label{revised}
In order to maximize \prob{S_n \leq t}, it seems that as much probability mass as possible should be on or near the
boundary $S_n=t$; the support condition from Lemma~\ref{support} illustrates this. Furthermore, the
fewer support points $\mu$ has, the more mass can contribute to the event $S_n=t$; an illustration
of this can be seen in the remark on page~\pageref{alternative}, where shrinking a continuous
portion of the distribution to one point doubles the contribution to \prob{S_n \leq t}.
Sometimes, however, putting some probability mass at $1$ may enable a redistribution of mass on lower support points that
results in an increase of \prob{S_n \leq t}. This (we think) is the intuitive explanation for the $\{0, t, 1\}$
solution.  It is also the reason we think that the number of support points required is no larger than three.

\begin{conjecture}
For any $n \geq 2$, $0<m<1$, and $0 \leq t < m n$, all distributions attaining the supremum $p_n(m,t)$ 
belong to the collections described below.
\end{conjecture}

\paragraph{Conjectured extremal binary solutions.}
A collection of at most $n$ distributions with two support points is identified below. They may not all satisfy all
of the Lagrange conditions.  
However, it is conjectured that if an extremal distribution is binary, it \emph{must} be one of these.

Suppose the support is $\{a, b\}$, with $0 \leq a <m<b\leq 1$. Lemma~\ref{support} implies $t-a \in \supp{S_{n-1}}$,
which means that $t-a = j \, b + (n-1-j) \, a$ for some integer $j=0, 1, \ldots, n-1$. From $\expec{X}=m$ follows that $\pi
:= \prob{X=b} = (m-a)/(b-a)= (m-a) \, j /(t - n \, a)$ (where $b-a = (t -n \, a)/j$ is used) and so
$\prob{S_n \leq t} =b(n, \pi, j)$.
Since $\pi$ is increasing in $a$ and $b(n,\pi,j)$ decreasing in $\pi$, one should minimize $a$.  
From $0 \leq a <m$ it follows that $m - (m\,n-t)/j < b \leq t/j$, so for $0 \leq j \leq t$ the constraint $b \leq
1$ becomes active as $a \downarrow 0$. Hence, for these $j$, the solution is $b_j=1$, $a_j=(t-j)/(n-j)$ (from
$a=(t-j\,b)/(n-j)$) and $\pi_j=1- (1-m) (n-j)/(n-t)$.
Considering that $m<b\leq 1$ implies $(t-j)/(n-j) \leq a < (t- j \, m)/(n-j)$, one sees that for $t<j<t/m$ one
should set $a_j=0$, $b_j=t/j$, and $\pi_j=j\,m/t$.

This results in a collection of at most $n$ potential extremal distributions, from which the best is
selected by comparing the values $b(n, \pi_j, j)$, for $j=0, 1, \ldots \ceil{t/m}-1$.

\paragraph{Conjectured extremal ternary solutions.}
It was shown on page~\pageref{zeroandone} that the supremum $p_n(m,t)$ is attained by a distribution with $0$ or $1$
in the support. The intuitive argument given above suggests that an extremal \emph{ternary} distribution will have \emph{both} 
$0$ and $1$ in the support. 
Using this as an assumption, a collection of (at most $\binom{n}{2}$) possible three point supports can be identified.
In order to precisely specify the distributions, the Lagrange linearity condition~(\ref{L2}) is needed as well.

Suppose the support is $\{0, a, 1\}$ with $0<a<1$. Lemma~\ref{support} implies $\{t-a, t\} \subset \supp{S_{n-1}}$,
whence integers $k\geq 1$ and $l\geq 0$ should exist, such that $k+l\leq n-1$ and $t=k \, a+ l$. Solving the last
equation for $a$, define $a_{k,l} = (t-l)/k$, which is between $0$ and $1$ if $0 \leq l <t < l+k \leq n-1$.
In contrast with the binary solutions above, the requirement that $\expec{X}=m$ is insufficient to fix the
probabilities and as an additional equation one should use the Lagrange linearity condition: \prob{S_{n-1} \leq t-x}
is linear for $x \in \{0, a, 1\}$. This results in
\beqs
(1-a) \, \prob{t-a < S_{n-1} \leq t} = a \, \prob{t-1 < S_{n-1} \leq t-a}.
\eeqs
Since $S_{n-1}$ is distributed as $a\,N_a+N_1$, where $(N_0,N_a,N_1)$ have a trinomial distribution with parameters
$n-1$, $p=\prob{X=1}$, $q=\prob{X=a}$ and $r=\prob{X=0}$, 
this last requirement is a polynomial equation in $p$, $q$, and $r$.
The requirements $a\,q+p=m$ and $p+q+r=1$ can be used to eliminate $q$ and~$r$, leaving a polynomial equations of
order $n-1$ in $p$.

\bibliographystyle{plain}
\bibliography{sr866}
\end{document}